\documentclass{article}
\usepackage{graphicx} %
\usepackage[dvipsnames]{xcolor}
\title{Sparse Pseudospectral Shattering}
\author{Rikhav Shah\footnote{Supported by NSF CCF-2420130}\\ UC Berkeley
       \and
       Nikhil Srivastava\footnote{Supported by NSF CCF-2420130}\\ UC Berkeley
       \and
        Edward Zeng\\ NYU}
\date{March 31, 2026}

\usepackage{graphicx} %
\usepackage{amsmath,amssymb,amsthm,amsfonts}
\usepackage{algpseudocode}
\usepackage{algorithm}
\usepackage{fullpage}
\usepackage{hyperref}
\usepackage{xcolor}
\usepackage{mleftright}
\usepackage{enumitem}
\usepackage{tikz}
\usetikzlibrary{shapes, arrows}

\newcommand{\net}{{\mathcal{N}}}
\newcommand{\C}{{\mathbb{C}}}
\newcommand{\R}{{\mathbb{R}}}
\newcommand{\Z}{{\mathbb{Z}}}

\renewcommand{\L}{{\mathcal{L}}}
\newcommand{\sph}[1]{{\mathbb{S}^{#1}}}

\newcommand{\eps}{\varepsilon}

\newcommand{\inr}[2]{{\left\langle#1,#2\right\rangle}}
\newcommand{\rinr}[2]{{\langle#1,#2\rangle}}

\newcommand{\st}{\,\textnormal{s.t.}\,}
\newcommand{\wrt}{\,\textnormal d}

\newcommand{\qand}{\quad\text{ and }\quad}

\newcommand{\abs}[1]{\mleft|#1\mright|}

\newcommand{\magn}[1]{\left\|#1\right\|}
\newcommand{\rmagn}[1]{\|#1\|}

\newcommand{\pare}[1]{\mleft(#1\mright)}

\newcommand{\sqbrac}[1]{{\left[{#1}\right]}}
\newcommand{\set}[1]{{\left\{{#1}\right\}}}

\newcommand{\bmat}[1]{\begin{bmatrix}#1\end{bmatrix}}

\newcommand{\alg}[1]{\textnormal{\texttt{#1}}}

\newcommand{\spliteq}[2]{\begin{equation}#1\begin{split}#2\end{split}\end{equation}}

\DeclareMathOperator*{\E}{\mathbb{E}}

\DeclareMathOperator*{\argmin}{arg\,min}

\DeclareMathOperator*{\poly}{poly}

\DeclareMathOperator*{\vol}{vol}

\DeclareMathOperator*{\supp}{supp}

\DeclareMathOperator*{\ball}{Ball}
\DeclareMathOperator*{\nnz}{nnz}

\newtheorem{theorem}{Theorem}[section]
\newtheorem{lemma}[theorem]{Lemma}
\newtheorem{proposition}[theorem]{Proposition}
\newtheorem{corollary}[theorem]{Corollary}
\newtheorem{fact}[theorem]{Fact}

\theoremstyle{definition}

\newtheorem{remark}[theorem]{Remark}

\algdef{SE}[REPEATN]{RepeatN}{End}[1]{\algorithmicrepeat\ #1 \textbf{times}}{\algorithmicend}

\usepackage{pifont}
\definecolor{darkblue}{rgb}{0.2, 0.2, .9}

\DeclareMathOperator{\Ber}{Bernoulli}

\newcommand{\bP}{\Pr}
\newcommand{\specr}{spr}

\begin{document}

\maketitle
\begin{abstract}
The eigenvalues and eigenvectors of nonnormal matrices can be unstable under
perturbations of their entries. This renders an obstacle to the analysis of
numerical algorithms for non-Hermitian eigenvalue problems. A recent technique
to handle this issue is pseudospectral shattering
\cite{b0}, showing that adding a random perturbation to any
matrix has a regularizing effect on the stability of the eigenvalues and
eigenvectors. Prior work has analyzed the regularizing effect of dense Gaussian
perturbations, where independent noise is added to every entry of a given
matrix \cite{b1,b0,b2,b3}. 

We show that the same effect can be achieved by adding a sparse random perturbation. In
particular, we show that given any $n\times n$ matrix $M$ of polynomially
bounded norm: (a) perturbing $O(n\log^6(n))$ random entries of $M$ by adding
i.i.d.~complex Gaussians yields $\log\kappa_V(A)=O(\poly\log(n))$ and $\log
(1/\eta(A))=O(\poly\log(n))$ with high probability; (b) perturbing
$O(n^{1+\alpha})$ random entries of $M$ for any constant $\alpha \in (0, 1/2]$ yields
$\log\kappa_V(A)=O_\alpha(\log(n))$ and $\log(1/\eta(A))=O_\alpha(\log(n))$ with
high probability. Here, $\kappa_V(A)$ denotes the condition number of the
eigenvectors of the perturbed matrix $A$ and $\eta(A)$ denotes its minimum
eigenvalue gap. 

A key mechanism of the proof is to reduce the study of $\kappa_V(A)$ to
control of the pseudospectral area and minimum eigenvalue gap of $A$, which are
further reduced to estimates on the least two singular values of shifts of $A$. We
obtain the required least singular value estimates via a streamlining of an
argument of Tao and Vu \cite{b4} specialized to the case of sparse complex
Gaussian perturbations.

Our main application is to an algorithm for solving $Ax=b$. We obtain a backward error guarantee for a modification of GMRES that uses $ O((\poly\log(n)+\log(1/\eps))(\nnz(A)+n\poly\log(n)+n\log(1/\eps)))$ operations under a natural spectral assumption,  namely when a certain pseudospectrum of $A$ is contained in a closed disk which is well-separated from the origin.

\end{abstract}

\section{Introduction}
A central question in numerical analysis is ``how do the eigenvalues and
eigenvectors of a matrix behave under perturbations of its entries?'' For normal
matrices, the eigenvalues are $1-$Lipschitz functions of the entries, and the
eigenvectors are stable under perturbations if the minimum eigenvalue gap is
large. This fact is essential to the rapid convergence and rigorous analysis of
algorithms for the Hermitian eigenvalue problem and its cousins.

For nonnormal matrices, two related difficulties appear: {\em
non-orthogonality} of the eigenvectors and {\em spectral instability}, i.e. high sensitivity of the
eigenvalues to perturbations of the matrix entries. Non-orthogonality slows down
the convergence of iterative algorithms (such as the power method) and spectral
instability makes it difficult to rigorously reason about convergence in the presence of
roundoff error. The main tool used to surmount these difficulties in recent years is smoothed analysis, i.e., adding a small
random perturbation to the input and solving the perturbed problem,
incurring a small backward error\footnote{i.e., the algorithm produces an exact
solution to a nearby problem, rather than producing an approximate solution to
the given problem, which is called forward error. }.
Specifically it was shown in \cite{b0} that adding small i.i.d. complex
Gaussian random variables to each entry of a matrix produces a matrix with well-conditioned
eigenvectors and a large eigenvalue gap, a phenomenon termed ``pseudospectral shattering.'' This 
was then generalized to other random variables in \cite{b1,b3,b5}, 
and is currently an essential mechanism in {all} of the known convergence results about diagonalizing arbitrary dense matrices in finite arithmetic \cite{b6,b0,b7,b8,b9}.

All existing works examine the setting where i.i.d. noise is added to every entry of a given matrix.
This paper asks if it possible to achieve pseudospectral shattering by adding noise to only a subset of entries, selected at random.
We provide a positive answer to this question in the regimes when the random perturbation has $O(n\log^6(n))$ and $O(n^{1+\alpha})$ nonzero entries.

Our results are phrased in terms of the sparsity $\rho=\rho(n)$ of the added noise. 
Given a matrix $M$, we consider perturbations of the form 
\begin{equation}\label{a0}A=M+N_g,\end{equation}
where the entries of $N_g$ are i.i.d.~copies of $\delta\cdot g$ where $\delta\sim\Ber(\rho)$ and $g\sim\mathcal N(0,1_\C)$.
As one might expect, our guarantee provides stronger regularization for larger $\rho$. We measure regularization in terms of the 
\textit{eigenvector condition number} $\kappa_V(A)$ and \textit{minimum eigenvalue gap} $\eta(A)$.
In the following definitions of these quantities, $A=VDV^{-1}$ is any diagonalization of $A$ and $\lambda_1(A),\ldots,\lambda_n(A)$ are the eigenvalues of $A$, counted with multiplicity:
\[
    \kappa_V(A) = \inf_{A = VDV^{-1}}\rmagn{V^{-1}}\magn V
\qand
    \eta(A)=\min_{i\neq j}\abs{\lambda_i(A)-\lambda_j(A)}.
\]
When $A$ is not diagonalizable, we state $\kappa_V(A)=\infty$.
Our main result is the following.
\begin{theorem}\label{a1}
For any $\rho$, $M \in \mathbb C^{n\times n}$ satisfying
\begin{equation} \label{a2}
    \frac{\left[\log(n)\log(2\|M\| + n)\right]^{3}}{n} \le \rho \le \frac{n^{1/2}}{n}
\end{equation}
we have
\begin{align}
    \Pr\pare{\kappa_V(M+N_g)\ge(2\magn M+n)^{15K}}
    &\le O\pare{n^{-K}} \label{a3} \\
    \Pr\pare{\eta(M+N_g)\le(2\magn M+n)^{-15K}}
    &\le O\pare{n^{-K}}, \label{a4} 
\end{align}
where $K=\lceil 3\log(n)/(2\log(n\rho))\rceil$.
\end{theorem}
\noindent Thus, we  obtain bounds of the form
\begin{equation}\label{a5}\log\kappa_V(M+N_g)\le O_\alpha(\log(n))\qand \log(1/\eta(M+N_g))\le O_\alpha(\log(n))\end{equation}
with high probability when $\magn M\le\poly(n)$ and $\rho(n)=n^\alpha/n$ for $\alpha \in (0, 1/2]$. In the sparser regime $\rho(n)=\log^6(n)/n$, we obtain the weaker bounds \spliteq{\label{a6}}{
    \log(\kappa_V(M+N_g))\le\poly(\log(n))\qand\log(1/\eta(M+N_g))\le\poly(\log(n)).
}

\begin{remark}[Small Perturbations]
One can derive a similar bound with smaller perturbations, i.e. $\sigma g$ in place of $g$, by noting that $\kappa_V(\sigma A)=\kappa_V(A)$ and $\eta(\sigma A)=\sigma\eta(A)$ for any scalar $\sigma>0$. We find it more convenient to work with $A=M/\sigma+N_g$ instead of $A=M+\sigma N_g$, so that the entries of $A-\E(A)$ have density bounded by 1 everywhere except at 0.
Using this observation gives the following Corollary of Theorem \ref{a1}.
\end{remark}
\begin{corollary}\label{a7}
Suppose $M\in\C^{n\times n}$ satisfies $\magn M\le1$. Let
\[
    \delta \ge \exp\left(-n^{1/7}\right), \quad \rho = \frac{\left[\log(n)\log(2/\delta+n)\right]^3}n, \quad
    \chi(n)=\left\lceil \frac{\log(n)}{(2\log\log n +2\log\log(n + 2/\delta)}\right\rceil.
\]
Then
\begin{align}
    \Pr\pare{\kappa_V(M+\delta N_g)\ge\pare{\frac2\delta+n}^{10\chi(n)}}
    &\le O\pare{n^{-\chi(n)}},\\
    \Pr\pare{\eta(M+\delta N_g)\le\delta\cdot\pare{\frac2\delta+n}^{-15\chi(n)}}
    &\le O\pare{n^{-\chi(n)}},
\end{align}
\end{corollary}

The lower bound on $\delta$ is due to the inequality  
\[
    \frac{\left[\log(n)\log(2\|M\| + n)\right]^{3}}{n} \le \frac{n^{1/2}}{n}
\]
in Theorem \ref{a1}, which is satisfied if $\|M\| \le \exp(n^{1/7})$.

\begin{remark}[Sharpness of the Bounds]
A coupon collector argument shows that, in the model \eqref{a0}, the sparsity condition $\rho(n)=\Omega(\log(n)/n)$ is required in order to obtain any high probability lower bound on $\eta(A)$. Theorem \ref{a1} requires a logarithmically denser perturbation than this to obtain \eqref{a6}.

On the other hand, the upper bound $\rho(n) \le 1/\sqrt{n}$ in Theorem \ref{a1} is not intrinsic. It arises only from our choice to specialize the least singular value estimates in \S\ref{a8} to the sparse regime. With a separate treatment of the sparse and dense cases in \S\ref{a8}, this could be relaxed to the full range $\rho(n)\le 1$. We do not pursue this here, in the interest of simplicity.

Finally, we have made no attempt to optimize the constants in our results.
\end{remark}

\subsection{Algorithmic Significance}
Our primary application is to an analysis of GMRES, an iterative algorithm for solving square, nonsymmetric, systems $Mx=b$ \cite{b10}. In contrast to other standard approaches such as using CG to solve $M^*Mx=M^*b$, GMRES can offer two advantages : (1) its convergence depends on the spectral properties of $M$ rather than $M^*M$, in some cases yielding faster convergence. (2) it can be implemented using only matrix-vector queries to $M$ as opposed to queries to both $M$ and $M^*$, which is important in certain PDE applications (see e.g. \cite{b11}).

When run for $k$ iterations, the algorithm outputs an approximate solution $x_k$ satisfying
\[
x_k=p_k(M)b,\quad p_k=\argmin_{\deg(p)\le k}\magn{b-Mp(M)b}
\]
where $p$ is a polynomial \cite{b10}. There are several bounds on the convergence of this algorithm, cataloged and compared by Embree in \cite{b12}. Roughly speaking, the bounds on the residual error have the form
\begin{equation}\label{a9}
\frac{\magn{b-Mx_k}}{\magn b}\le C_{M,\Omega}\beta_{\Omega, k}
\end{equation}
where $\Omega\subset\C$ is a set containing the spectrum of $M$, $C_{M,\Omega}$ is a constant depending on $M$ and the choice of $\Omega$, and $\beta_{\Omega,k}$ is related to the polynomial approximation problem on $\Omega$, namely,
\[\beta_{\Omega,k}:=\min_{\deg(p)\le k}\sup_{z\in\Omega}\abs{1-zp(z)}.\]
When $\Omega$ is large, the polynomial approximation problem becomes more difficult, so $\beta_{\Omega,k}$ increases, but sometimes with the benefit of decreasing $C_{M,\Omega}$. At one extreme end where $\Omega$ is as small as possible, we have the original bound of Saad and Schultz \cite{b10} which takes $\Omega=\Lambda(M)$ to be the spectrum of $M$ and shows $C_{M,\Lambda(M)}=\kappa_V(M)$ satisfies \eqref{a9}. At the other extreme end where $C_{M,\Omega}$ is as small as possible, \eqref{a9} is satisfied when $\Omega=\set{x^*Mx:\magn x=1}$ is the field of values or numerical range of $M$ and $C_{M,\Omega}=1+\sqrt2$ (the bound of $C_{M,\Omega}=O(1)$ is implied by the main result of Crouzeix \cite{b13}, and $1+\sqrt2$ by a follow-up work \cite{b14}), c.f. \cite{b15,b12}.

This work proposes a one-line modification of the GMRES algorithm: add a sparse random perturbation to the initial input before running the standard GMRES algorithm. The resulting bound proved by this work can be seen as an interpolation between the extreme ends of $C_{M,\Omega}$ versus $\Omega$, attempting to capture the attractive features of both. Our guarantee is slightly weaker than \eqref{a9}. Note that \eqref{a9} can be interpreted as a backward error guarantee. That is, it says that the output $x_k$ solves exactly a perturbed input
\[Mx_k=(b+\delta_b)\]
for the error bound $\magn{\delta_b}\le \magn bC_{M,\Omega}\beta_{\Omega,k}$.
Our modification of GMRES solves exactly a perturbed system
\[(M+\Delta_M)x_k=(b+\delta_b)\]
for the error bound
\[\frac{\magn{\Delta_M}}{\magn M}\le\eps,\quad\frac{\magn{\Delta_b}}{\magn b}\le C_{M,\Omega}\beta_{\Omega,k}.\]

In Theorem \ref{a10}, we show that this approach achieves a value of $C_{M,\Omega}$ which is independent of $M$, with $\Omega=\Lambda_{\eps}(M)$, whereas previous results with this choice of $\Omega$ had at least an implicit dependence of $C_{M,\Omega}$ on $M$. Notably, Trefethen and Embree \cite{b16,b17}, c.f. \cite{b12} show \eqref{a9} where $C_{M,\Omega}$ is the perimeter of $\Lambda_\eps(M)$ divided by $\eps$. By the Bauer-Fike theorem, this is upper bounded by $O(\kappa_V(A))$, with a matching close lower bound for very small $\eps$. %

The algorithm is best suited to cases where $\beta_{\Omega,k}$ is an exponentially decreasing quantity in $k$, which it is in our case under the natural assumption that $\Lambda_\epsilon(M)$ is contained in a disk which is well-separated from the origin (see Theorem \ref{a10} and the discussion following it). In these settings, the number of iterations required to achieve $\delta$ backward error is just $O\pare{\log(C_{M,\Omega}/\delta)}$, which after our perturbation will be just $O_\beta(\log^2(n)+\log(1/\delta))$. Without a perturbation, amortized cost of each iteration of GMRES is $O(nk+\nnz(M))$, where $\nnz(M)$ is the number of nonzero entries of $M$. Since our perturbation is sparse, the first step of our algorithm increases $\nnz(M)$ by at most $O(n\log^6(n))$.

A second application of Theorem \ref{a1} is to  \textit{derandomization}. Earlier in the introduction, we identified diagonalization algorithms as a motivation for pseudospectral shattering.
Those earlier efforts identified a polynomial dependence on $\log(\kappa_V(A)/\eta(A))$ in the complexity of a fast diagonalization algorithm. This quantity is at most roughly $O(\log(n))$ for dense random perturbations of the input, and at most roughly $O(\log^2(n))$ for sparse perturbations by Theorem \ref{a1}. Thus, for only additional factor of $\log(n)$ overhead, we reduce the number of random bits required by the algorithm from roughly $O(n^2\log n)$ to $O(n\log^6 n)$ (in the regime when the desired backward error is $1/\poly(n)$), since $\log(n)$ bits are used for each nonzero entry of the perturbation.
In the context of numerical analysis, where replicability is desirable, it means that the random bits used by the algorithm can be stored in nearly linear rather than quadratic space.

\subsection{Related Work} 
\paragraph{Random Matrix Theory.} Several works have studied the distributions of $\sigma_n(\cdot),\kappa_V(\cdot),$ and $\eta(\cdot)$ for $n\times n$ random matrices $A$ where $A-\E(A)$ has i.i.d. entries. The table below lists the relevant works and which quantities they consider, under the normalization $\magn{\E(A)}\le \poly(n)$.
These works all take the entries of $N_x=A-\E(A)$ to be of the form $\delta\cdot x$ where $\delta\sim\Ber(\rho)$ is a Bernoulli random variable with parameter $\rho=\rho(n)$, and $x$ is a random variable which varies from work to work. 
\newcommand{\xmark}{}
\begin{table}[h]
\centering
\begin{tabular}{l|c|c|c|c|c|c}
 & $\kappa_V(\cdot)$ & $\eta(\cdot)$ & $\sigma_n(\cdot)$ & $\E(A)$ & $x$ & $\rho(n)$ \\
\hline
\hline
\cite{b4} & \xmark & \xmark & \checkmark & any & finite moment conditions & $n^\alpha/n$ \\
\cite{b18} & \xmark & \xmark & \checkmark & 0 & finite moment conditions & $\log^c(n)/n$ \\
\cite{b19} & \xmark & \xmark & \checkmark & $zI$ & finite moment conditions & $\omega(1/n)$ \\
\cite{b20} & \xmark & \checkmark & \checkmark & 0 & finite moment conditions & $\Theta(1)$ \\
\cite{b2} & \checkmark & \xmark & \checkmark & any & complex Gaussian & 1 \\
\cite{b0} & \checkmark & \checkmark & \checkmark & any & complex Gaussian & 1 \\
\cite{b1,b3} & \checkmark & \checkmark & \checkmark & real & bounded real density & 1 \\
\cite{b5} & \checkmark & \checkmark & \checkmark & real/complex & real/complex Gaussian & 1 \\
\hline
This work & \checkmark & \checkmark & \checkmark & any & complex Gaussian & $\Omega(\log^2(n)/n)$ \\
\hline
\end{tabular}
\end{table}
From the table, one can see that this work provides the first bounds on $\kappa_V(\cdot)$ and $\eta(\cdot)$ for $\mathbb E(A)\neq 0$ and $\rho(n)$ less than $1$. In particular, it is the first work on those quantities when $\mathbb E(A) \neq 0$ and the entries of $A-\E(A)$ do not have absolutely continuous distributions on either $\R$ or $\C$.\\

\noindent {\em Bounds on $\kappa_V(\cdot)$ and $\eta(\cdot)$.} The listed prior works for $\kappa_V(\cdot)$ and $\eta(\cdot)$ obtain bounds of the form $\kappa_V(M+N_x)\le\poly(n)$ and $\eta(M+N_x)\ge\poly(n)^{-1}$ with high probability for $\|M\|\le \poly(n)$. The series of works \cite{b21}, \cite{b22}, culminating in \cite{b5} listed in the table, obtains bounds for real and complex Gaussian perturbations with nearly optimal dependence on $n$.
All of these works bound $\kappa_V(A)$ by obtaining good control over the size of the {\em $\eps$-pseudospectrum} of $A$, defined as
\[\Lambda_\eps(A)=\set{z\in\C:\sigma_n(z-A)\le\eps}.\]
It is well-known that the set $\Lambda_\eps(A)$ always contains disks of radius $\eps$ around each eigenvalue of $A$, and that equality is achieved if and only if $\kappa_V(A)=1$ \cite{b17}. A generalization of this fact shows that the limit of $\eps^{-2}\vol\Lambda_\eps(A)$ as $\eps\to0$ determines $\kappa_V(A)$ up to a $\poly(n)$ factor (see \cite[Corollary 1.5]{b23} and \cite[Lemma 3.2]{b2}). In order for this limit to be finite, one must show the that $\E \vol\Lambda_\eps(A)=O_n(\eps^2)$ as $\eps\rightarrow 0$, which follows from least singular value estimates on shifts of $A$ with the correct exponent, namely:
\begin{equation}\label{a11} \Pr\pare{\sigma_n(z-A)\le \eps}=O_n(\eps^2)\quad\forall \eps>0, z\in\C.\end{equation}
\cite{b2} prove \eqref{a11} in the dense complex Gaussian case, yielding a bound on $\kappa_V(A)$; \cite{b0} then uses this to control $\eta(A)$ as well.
In the more challenging setting of real Gaussian perturbations, the bound \eqref{a11} does not hold for $z\in\R$. The subsequent works \cite{b1,b3,b5} overcome this obstacle by first controlling $\eta(A)$, which then allows them to use some alternate arguments which handle $z\in\R$ and $z\in\C\setminus\R$ separately. One such argument is the bootstrapping scheme of \cite{b3}, which shows that an adequately good lowerbound on $\eta(A)$ and an adequately good upperbound on $\E\vol\Lambda_\eps(A)$ for specific scales $\eps>0$ (not tending to zero) together imply a bound on $\kappa_V(A)$. While their approach was originally designed to handle the specific difficulties arising in the case of real Gaussian perturbations, we are able to repurpose it  sparse complex Gaussian perturbations (see Section \ref{a12} for details). 

We remark that obtaining the exponent of $2$ in \eqref{a11} for complex $z$ (or something close to it) was a sticking point in all previous works. A key technical contribution of this work, described in detail in Section \ref{a12}, is that we show how to obtain bounds on $\kappa_V(A)$ from much weaker least singular value bounds: rather than \eqref{a11}, it suffices to have a bound of type \eqref{a13} satisfying just \eqref{a14}. Our approach thus reveals that the main bottleneck in proving bounds on $\kappa_V(A)$ is proving bounds on $\eta(A)$. We view this realization as a conceptual contribution of this work.  \\

\noindent {\em Bounds on $\sigma_{n}(\cdot)$.}The result in the first row of the above table was provided by Tao and Vu in the context of proving the circular law \cite[Theorem 2.9]{b4}. They showed $\sigma_n(M+N)\ge n^{-b}$ with probability at least $1-n^{-a}$ where $b$ depends on $a$ and $\alpha$.
This estimate was extended to lower sparsity $\rho(n)=\log^c(n)/n$ with the restriction that $\E(A)=0$ by Basak and Rudelson \cite{b18}, and pushed even further down to any $\rho(n)=\omega(1/n)$ with the milder restriction $\E(A)=zI$ by 
Rudelson and Tikhomirov \cite{b19} (though in that case, $b$ must be replaced by $b\log^2(n)$).

As explained in Section \ref{a12}, bounds on the least singular
value play a key role in our results, and we use many ideas from \cite{b4}.
Our technical innovation here is that by restricting attention to complex
Gaussian $x$, we are able to avoid the use of additive combinatorics and prove
stronger bounds than previous works, which are needed in this context. In particular, we obtain an improved
dependence $b=a/2+O(1)$ in the Tao-Vu bound.

\paragraph{Numerical Analysis.} The idea of using a small random perturbation to
improve the performance of numerical algorithms (specifically, Gaussian
elimination) originated in the influential work \cite{b24}, which analyzed the least
singular value under a dense Gaussian perturbation. The first paper to propose
the use of a (not necessarily random) perturbation to tame spectral instability in the context of diagonalization algorithms was \cite{b25}.

\cite{b6} gave the first numerically stable algorithm for
diagonalizing a dense non-Hermitian matrix; a key step there was to add a dense
complex Gaussian perturbation to the input, and show that it yields a 
well-conditioned diagonalization problem in a certain geometric sense. The more recent works
\cite{b0,b8,b9} used pseudospectral shattering in their analyses of nearly matrix multiplication time and $O(n^3)$ time algorithms for computing the eigendecomposition and Schur form of a matrix. This technique was extended to computing invariant subspaces of Hermitian pencils by \cite{b26}, to diagonalizing arbitrary matrix pencils by \cite{b27}, and to computing the Schur form of arbitrary matrix pencils by \cite{b28}.

\newcommand{\minihead}[1]{\vspace{3mm}\noindent\emph{#1:}}
\subsection{Technical Overview}\label{a12}
We now describe the proof of Theorem \ref{a1}.
The proof consists of three steps. In the first two steps, $A$ can be any random matrix.
In the third step, we need our particular model of sparse perturbations, $A=M+N_g$. When we omit the subscript we mean $N=N_g$. We dedicate a section to each step.
The bulk of the technical work is in the first and third steps; the second step can be viewed as the glue which combines them to produce the final theorem. The logical structure of this paper is captured by the following diagram.\\

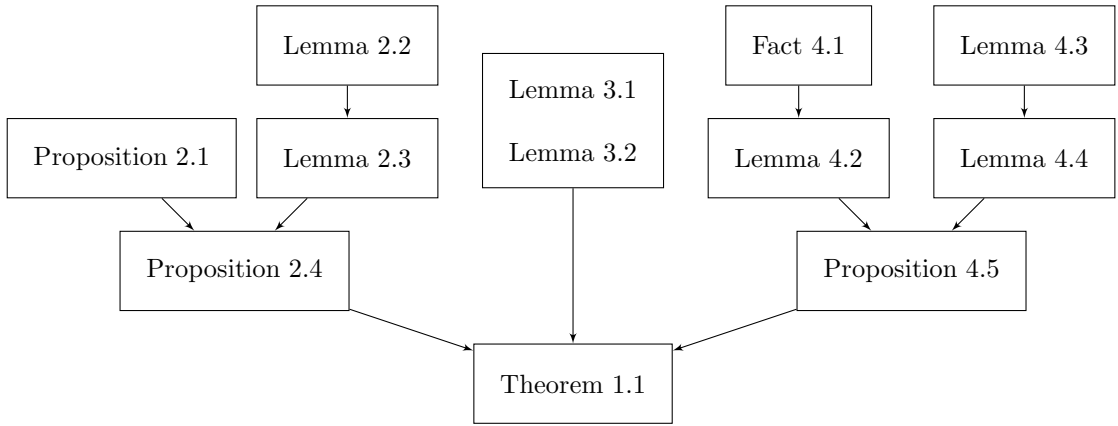
\tikzstyle{terminator} = [rectangle, draw, text centered, rounded corners, minimum height=2em]
\tikzstyle{process} = [rectangle, draw, text centered, minimum height=2em]
\tikzstyle{decision} = [diamond, draw, text centered, minimum height=2em]
\tikzstyle{data}=[trapezium, draw, text centered, trapezium left angle=60, trapezium right angle=120, minimum height=2em]
\tikzstyle{connector} = [draw, -latex']

\tikzstyle{thm}=[rectangle, draw, text centered,inner sep=1em, minimum height=3em]

\begin{tikzpicture}
\node [thm] at (-6,2.5) (p21) {Proposition \ref{a15}};
\node [thm] at (-3,4) (l22) {Lemma \ref{a16}};
\node [thm] at (-3,2.5) (l23) {Lemma \ref{a17}};
\node [thm] at (-4.5,1) (p24) {Proposition \ref{a18}};
\path [connector] (p21) -- (p24);
\path [connector] (l22) -- (l23);
\path [connector] (l23) -- (p24);

\node [thm] at (3,4) (f41) {Fact \ref{a19}};
\node [thm] at (3,2.5) (l42) {Lemma \ref{a20}};
\node [thm] at (6,4) (l43) {Lemma \ref{a21}};
\node [thm] at (6,2.5) (l44) {Lemma \ref{a22}};
\node [thm] at (4.5,1) (p45) {Proposition \ref{a23}};
\path [connector] (f41) -- (l42);
\path [connector] (l42) -- (p45);
\path [connector] (l43) -- (l44);
\path [connector] (l44) -- (p45);

\node [thm, align=center] at (0,3) (l3x) {Lemma \ref{a24} \\ \\ Lemma \ref{a25}};

\node [thm, align=center] at (0,-.5) (main) {Theorem \ref{a1}};

\path [connector] (p24) -- (main);
\path [connector] (l3x) -- (main);
\path [connector] (p45) -- (main);
\end{tikzpicture}

\minihead{Step 1. Reduction to pseudospectral area and eigenvalue gap (Section \ref{a26})}
We first show (Proposition \ref{a18}) that an upper tail bound on $\Pr(\kappa_V(A)\ge \tau)$ may be derived from a lower tail bound on $\eta(A)$ {\em and} a bound of type 
\begin{equation}\label{a27}\E\vol\Lambda_\eps(A)\le\poly(n)\cdot\eps^c+
\exp\left(-\Omega(\log^2(n)/\log(n\rho))\right),\end{equation}
for any constant $c>0$, where $\eps$ is a carefully chosen spectral scale depending $c$ and on the desired deviation $\tau$. The proof is an adaptation of the bootstrapping scheme of \cite{b3} and additionally invokes a deterministic exponential upperbound on $\kappa_V(A)$ in terms of $\eta(A)$ (Proposition \ref{a15}).
In particular, the allowance of the additive term is essential in the sparse Gaussian case as there is some nontrivial probability that one or more rows and columns of the matrix end up with no noise added to them.

\minihead{Step 2. Reduction to bottom two singular values (Section \ref{a28})}
The previous step shows we need control over $\E\vol\Lambda_\eps(A)$ and $\eta(A)$.
As mentioned above, $\E\vol\Lambda_\eps(A)$ reduces to lower tail estimates for the least singular value $\sigma_n(z-A)$ for $z\in\C$ (see Lemma \ref{a24}).
Following \cite{b1,b3}, a lower tail bound on $\eta(A)$ can be reduced to lower tail bounds on the {\em
bottom two} singular values $\sigma_n(z-A),\sigma_{n-1}(z-A)$ via Weyl majorization and a standard net argument in the complex plane (see Lemma \ref{a25}).

The strength of a tail bound on $\sigma_{n-m}$ can be characterized in terms of the \textit{power} $c_m$ of $\eps$ on the right-hand side of a bound of the form
\spliteq{\label{a13}}{
\Pr\pare{\sigma_{n-m}(z-A)\le\eps}\le\poly(n)\eps^{c_m}+\beta(n),
}
where $\beta(n)$ is an additive error depending on $n$. In order to obtain a nontrivial bound on $\eta(A)$, our argument requires \eqref{a13} to hold for $m=0,1$ with
\begin{equation}\label{a14} \frac1{c_0}+\frac1{c_1}<1.\end{equation}

\minihead{Step 3. Control on bottom two singular values (Section \ref{a8})}
We show the required control \eqref{a13},\eqref{a14} over the bottom two singular values holds with room to spare. Specifically, we show a bound of the form \eqref{a13} holds for $c_0=2$ and $c_1=4$ with $\beta(n)=\exp(-\Omega(\rho n))$. (in fact, we show it holds for $c_m=2m+2$ for any constant $m$, see Proposition \ref{a23}).
The argument is based on an $\eps$-net construction following the compressible/incompressible or rich/poor decomposition in \cite{b4}.
That work considers lower tail bounds of $\sigma_n(M+N_x)$ where $x$ is a general subgaussian random variable and the sparsity parameter $\rho(n)=n^{\alpha-1}$, $\alpha>0$ is a fixed polynomial in $n$. They show for every $a>0$ there exists a $b>0$ such that
\[\Pr\pare{ \sigma_n(A)\le n^{-a} } \le n^{-b}.\]
By tracing their argument, one can show there is a linear relationship between $a$ and $b$ so that their bound more closely resembles the form \eqref{a13} for an unspecified tiny constant $c_0\le b/a$ and $\eps=1/\poly(n)$. This gives enough control over $\E\vol\Lambda_\eps(A)$ to satisfy \eqref{a27} in Step 2, but does not satisfy \eqref{a14} and is therefore not enough to bound $\eta(A)$ (see Remark \ref{a29}).

We resolve this issue by specializing to the complex Gaussian case $x=g$ and achieve three advantages over \cite{b4} in this setting.
Firstly, our argument applies to every $m$ (not just $m=0$), and we are able to obtain $c_0$ and $c_1$ satisfying \eqref{a14} (as a bonus, our proof yields the optimal power of $c_0=2$ for $m=0$).
Secondly, we obtain tail bounds at scales smaller than inverse polynomial in $n$.
Lastly, we are able to push the sparsity parameter down to $O(\log^6(n)/n)$. Furthermore, because $g$ has a continuous density, we avoid the additive combinatorics required by \cite{b4}, resulting in much simpler proofs.

In summary, while our proof uses adaptations of arguments from \cite{b3,b4,b1} as ingredients, it combines them in a new and more efficient way which clarifies their limits and allows us to handle the sparse setting.
\begin{remark}[Subgaussian and non-Gaussian perturbations]\label{a29}
Much of the proof of Theorem \ref{a1} goes through for \textit{any} random matrix model, not just $A=M+N_g$. When one replaces complex Gaussian $g$ with a subgaussian random variable, part of the proof of Theorem \ref{a1} (Lemma \ref{a18} in particular) can be combined with the existing lower tail bound on $\sigma_n(A)$ from Theorem 2.1 of \cite{b4} to produce a result of the form
\[
    \Pr\big(\kappa_V(M+N) \le \poly(n)\big) \ge \Pr\pare{\eta(M+N) \ge \poly(n)^{-1}} - \poly(n)^{-1}
\]
when $\log(n)/\log(n\rho)$ is a constant and $\magn M\le\poly(n)$. The reason we are not able to fully generalize Theorem \ref{a1} to the subgaussian setting is that we do not know how to prove strong enough tail lower tail estimates on $\eta(M+N)$.

All of our results hold if complex Gaussians are replaced by any sufficiently concentrated absolutely continuous distribution on $\C$, with bounds depending on the $\ell_\infty$ norm of the density and the decay of the tails. We have chosen to write the paper in terms of complex Gaussians for ease of exposition.
\end{remark}

\paragraph{Notation.}
$\wedge$ denotes logical ``and'' (intersection of events); $\vee$ denotes logical ``or'' (union of events); $\neg$ denotes logical negation (complement of event); $\mathbb S^{n-1}$ denotes the complex sphere $\set{z\in\C^n:\magn z=1}$; the matrix norms $\magn{\cdot}=\magn{\cdot}_2,\magn\cdot_{\infty}$ are the operator norms. $spr(\cdot)$ is the spectral radius.
Matrices are $n \times n$ unless otherwise specified.

\section{Reduction to eigenvalue spacing and pseudospectral area}
\label{a26}
Since a matrix with distinct eigenvalues must be diagonalizable, $\eta(A) > 0$ implies $\kappa_V(A)<\infty$. We begin with a quantitative version of this statement: $\kappa_V(A)$ is at most exponential in $1/\eta(A)$. 
\begin{proposition}[Exponential $\kappa_V$ bound]\label{a15}
Let $\magn A\le1$ and $\eta=\eta(A)$. Then\[\kappa_V(A)\le n\cdot 2^n\cdot \eta^{1-n}.\]
\end{proposition}
\begin{proof}
The case $\eta = 0$ is trivial so suppose $\eta > 0$. Let $A=VDV^{-1}$ for $V$ with unit column vectors. Number the eigenvalues $\lambda_1,\ldots,\lambda_n$ so that $D_{ii}=\lambda_i$.
Let\[d_J=\inf\set{\magn{Vx}:\magn x=1,\,\supp(x)\subset J}\]
so that $d_{[n]}=\sigma_n(V)$ and $\kappa_V(A)\le\sqrt n/d_{[n]}$. Let $x^J$ be a unit vector such that $d_J=\magn{Vx^J}$.
We argue via induction that
\[
    d_J\ge\frac{\eta^{\abs J-1}}{2^{\abs J}\sqrt{\abs J}}.
\]
The base case when $J$ is a singleton follows because the columns of $V$ are unit vectors.
Denote $i_J(x)=\argmin_{j\in J}\abs{x_j}$ . Then since $\magn A\le1$ we have
\spliteq{}{
d_J=\magn{Vx^J}\ge\magn{AVx^J}=\magn{VDx^J}
  &\ge\magn{V(D-\lambda_{i_J(x^J)}I)x^J}-\abs{\lambda_{i_J(x^J)}}\magn{Vx^J}
\\&\ge\magn{V(D-\lambda_{i_J(x^J)}I)x^J}-d_J.
}
Consider the vector $(D-\lambda_{i_J(x^J)}I)x^J$. It is supported on $J-\set{i_J(x)}$, and each absolute entry is at least $\eta$ times the corresponding entry of $x^J$ in absolute value. Thus its magnitude is at least
\[\eta\sqrt{1-\min_{j\in J}\abs{x_j}^2}\ge\eta\sqrt{1-\frac{1}{\abs J}}\] Therefore,
\[2d_J\ge\eta\sqrt{1-\frac1{\abs J}}\cdot d_{J-i_J(x^J)} = \eta\frac{\sqrt{\abs J}}{\sqrt{\abs{J}-1}}\cdot d_{J-i_J(x^J)}\]
which establishes the inductive step.
\end{proof}
For a more refined result, we define the \textit{eigenvalue condition numbers}.
Suppose $A$ has distinct eigenvalues with spectral expansion $A=\sum_j\lambda_jv_jw_j^*$. The eigenvalue condition number of $\lambda_j$ is defined as
\[
    \kappa(\lambda_j)=\magn{v_jw_j^*}.
\]
Bounds on the eigenvalue condition numbers imply bounds on $\kappa_V(A)$ as
\begin{equation} \label{a30}
    \max_j\kappa(\lambda_j) \le \kappa_V(A)\le\sqrt{n\sum_{j=1}^n\kappa(\lambda_j)^2}\le n\max_j\kappa(\lambda_j).
\end{equation}
It is well-known that the eigenvalue condition numbers of a matrix specify the rate at which the $\eps$-pseudospectrum shrinks in the $\eps \downarrow 0$ limit. The following lemma is a non-asymptotic version of this result.
\begin{lemma}[Disks inside pseudospectrum]
    \label{a16}
    Suppose $A$ has distinct eigenvalues $\lambda_1,\cdots,\lambda_n$ ordered such that $\kappa(\lambda_1)=\max_j\kappa(\lambda_j)$.
    Then
    \[
        \Lambda_\eps(A) \supset
        \ball\pare{\lambda_1,\frac12\min\pare{\frac{\eta(A)}n,\kappa(\lambda_1)\eps}}.
    \]
\end{lemma}
\begin{proof}
Let $\sum_j\lambda_jv_jw_j^*$ be the spectral expansion of $A$.
Let $r < \eta = \eta(A)$ and $z\in\ball(\lambda_1,r)$. Then by the triangle inequality,
\[
\magn{(zI-A)^{-1}}=\magn{\sum_j\frac{v_jw_j^*}{z-\lambda_j}}
\ge\frac{\kappa(\lambda_1)}r-\sum_{j\ge 2}\frac{\kappa(\lambda_j)}{\eta-r}
\ge\kappa(\lambda_1)\cdot\pare{\frac1r-\frac{n-1}{\eta-r}}=:f(r).
\]
Observe
\[
r(\eta-r)\pare{f(r)-\frac{1}{\eps}}
=
\kappa(\lambda_1)\eta-r\cdot\pare{n\kappa(\lambda_1)+\frac\eta\eps} + \frac{r^2}{\eps}
\ge\kappa(\lambda_1)\eta-r\cdot\pare{n\kappa(\lambda_1)+\frac\eta\eps}.
\]
The right-hand side is nonnegative when
\[
r
\le\frac{\kappa(\lambda_1)\eta}{\kappa(\lambda_1)n + \eta\cdot\frac1\eps}
  =\pare{\pare{\frac\eta n}^{-1} +\pare{\kappa(\lambda_1)\eps}^{-1}}^{-1}
\]
On the other hand, nonnegativity implies $f(r)\ge 1/\eps$, which implies $z\in\Lambda_\eps(A)$.
Therefore, if \[r\le\frac12\min\pare{\frac\eta n,{\kappa(\lambda_1)\eps}}\] then $\ball(\lambda_1, r) \subset \Lambda_\eps(A)$. Lemma \ref{a16} follows.
\end{proof}

We will use Lemma \ref{a16} to apply a bootstrapping argument to bound $\kappa_V$.
The following lemma is a single tug of the bootstrap, and closely follows \cite{b3}.
\begin{lemma}[Single bootstrap tug]\label{a17}
Let $A$ be any random matrix and let $F(\tau)$ denote the event that $A$ has distinct eigenvalues $\lambda_1, \dots \lambda_n$ with
\[
    \max_j\kappa(\lambda_j)^2\le\tau.
\]
Fix any $\tau_1<\tau_2$ and $\eps,\eta_0>0$.
If $\eta_0\ge n\eps\sqrt{\tau_2}$
then
\[
\Pr\pare{F(\tau_1)\wedge\eta(A)\ge\eta_0}
\ge
\Pr\pare{F(\tau_2)\wedge\eta(A)\ge\eta_0}
-\frac{4}{\pi} \cdot \frac{\E\vol\Lambda_\eps(A)}{\eps^2\tau_1}.
\]
\end{lemma}
\begin{proof}Without loss of generality take $\kappa(\lambda_1)=\max_j\kappa(\lambda_j)$.
    By Markov's inequality,
    \[
\Pr\pare{\neg F(\tau_1)\,|\,F(\tau_2)\wedge\eta(A)\ge\eta_0}
\le\frac{\E\pare{\kappa(\lambda_1)^2\,|\,F(\tau_2)\wedge\eta(A)\ge\eta_0}}{\tau_1}.
\]
Conditioned on $F(\tau_2)$, we have $\eta_0\ge n\eps\sqrt{\tau_2}\ge n\eps\kappa(\lambda_1)$, so by Lemma \ref{a16} we have
\[\vol\Lambda_\eps(A)\ge\frac\pi4\kappa(\lambda_1)^2\eps^2.\]
Rearranging gives
\spliteq{}{
\Pr\pare{\neg F(\tau_1)\,|\,F(\tau_2)\wedge\eta(A)\ge\eta_0}
  &\le\frac{\frac4\pi\E\pare{\vol\Lambda_\eps(A)\,|\,F(\tau_2)\wedge\eta(A)\ge\eta_0}}{\eps^2\tau_1}
\\&\le\frac{\frac4\pi\E\vol\Lambda_\eps(A)}{\eps^2\tau_1\cdot\Pr\pare{F(\tau_2)\wedge\eta(A)\ge\eta_0}}.
}
On the other hand,
\[
\Pr\pare{F(\tau_1)\,|\,F(\tau_2)\wedge\eta(A)\ge\eta_0}
\le
\frac
{\Pr\pare{F(\tau_1)\wedge\eta(A)\ge\eta_0}}
{\Pr\pare{F(\tau_2)\wedge\eta(A)\ge\eta_0}}.
\]
These two expressions sum to at least 1. Rearranging gives the result.
\end{proof}

\begin{proposition}[Complete bootstrap]\label{a18}
Let $A$ be any random matrix. Fix any $\tau_0 > 0$ and $R > \eta_0 > 0$. Let $\eps = \eta_0/\pare{n\tau_0^{1/2}}$.
Then
\spliteq{}{
\Pr\pare{\max_j\kappa(\lambda_j)^2\le\tau_0}
\ge
\Pr\pare{\eta(A)\ge\eta_0}
-\sqbrac{3n\lg\frac{3R}{\eta_0} + 3\lg\frac{2\eta_0}{\tau_0}}
\cdot\frac{\E\vol\Lambda_\eps(A)}{\tau_0\cdot\eps^2}
-\Pr\pare{\magn{A}>R}.}
\end{proposition}
\begin{proof}
    The idea is to apply Lemma \ref{a17} with a sequence of thresholds $\tau$.
    Set
    \[\tau_j=\tau_0\cdot 2^j\qand\eps_j=\frac{\eta_0}{n{\tau_j}^{1/2}}.\]
    As before, let $F(\tau)$ be the event $\max_j\kappa(\lambda_j)^2\le\tau$. Let $G$ be the event $\eta(A)\ge\eta_0$.
Then Lemma \ref{a17} implies
\spliteq{}{
\Pr(F(\tau_j)\wedge G)
  &\ge\Pr(F(\tau_{j+1})\wedge G)-\frac{4}{\pi} \cdot \frac{\tau_{j+1}}{\tau_{j}}\cdot\frac{n^2}{\eta_0^2}\cdot\E\vol\Lambda_{\eps_{j+1}}(A).
}
Note $\tau_{j+1}/\tau_j=2$. So by recursively applying Lemma \ref{a17} for $j=1,\cdots,K$, we obtain
\spliteq{}{\label{a31}
\Pr\pare{F(\tau_0)\wedge G}
  &\ge\Pr\pare{F(\tau_K)\wedge G}-\frac{8}{\pi} \cdot \frac{n^2}{\eta_0^2}\E\sum_{j=1}^{K}\vol\Lambda_{\eps_j}(A).
\\&\ge\Pr\pare{F(\tau_K)\wedge G}-\frac{8}{\pi} \cdot \frac{n^2}{\eta_0^2}\cdot K\cdot\E\vol\Lambda_{\eps_0}(A).}
Set
\[
    K = \left\lceil{n\lg\frac{3R}{\eta_0}+\lg\frac{\eta_0}{\tau_0}}\right\rceil \le n\lg\frac{3R}{\eta_0}+\lg\frac{2\eta_0}{\tau_0}.
\]
Then $\tau_K\ge \eta_0(3R/\eta_0)^n$ and so Proposition \ref{a15} combined with (\ref{a30}) implies $G\wedge\magn{A}\le R\implies F(\tau_K)$.
Hence
\spliteq{}{\label{a32}
\Pr(F(\tau_K)\wedge G)
  &\ge\Pr(F(\tau_K)\wedge G\wedge\magn{A}\le R)
\\&=\Pr(F(\tau_K)\,|\,G\wedge\magn{A}\le R)\cdot\Pr(G\wedge\magn{A}\le R)
\\&=\Pr(G\wedge\magn{A}\le R)
\\&\ge\Pr(G)-\Pr(\magn{A}>R).
}
The result follows from combining (\ref{a31}), (\ref{a32})
and writing $\eta_0$ in terms of $\eps_0,\tau_0$.
\end{proof}

\section{Reduction to singular values}
\label{a28}
In this section we reduce bounds on $\E\vol\Lambda_\eps(A)$ and $\eta(A)$ to estimates on small singular values.
\begin{lemma}[Pseudospectral area to smallest singular value]\label{a24}
Fix $M\in\C^{n\times n}$ and $R\ge\max(2\magn{M},n/2)$. Let $N$ be a random matrix with i.i.d.~symmetric subgaussian entries. Then
\[\E\vol\Lambda_\eps(M+N)
\le\pi (R + \eps)^2\sup_{\magn{M'}\le1.5R+\eps}\Pr\pare{\sigma_n(M'+N)\le\eps}
+\exp\pare{-\Omega(R)}.
\]
\end{lemma}
\begin{proof}
By Fubini,
    \spliteq{}{
    (\star):=\E\vol\Lambda_\eps(M+N)
      &=\E\int_{\C}\begin{cases}1&z\in\Lambda_\eps(M+N)\\0&\text{o.w.}\end{cases}\wrt z
    \\&=\int_{\C}\Pr\pare{\sigma_n(A-zI)\le\eps}\wrt z.}
Note $\sigma_n(M+N-zI)\ge\abs z-\magn{M+N}$. Thus
\spliteq{}{
(\star)
  &\le \pi (R + \eps)^2\sup_{\abs z\le R+\eps}\Pr\pare{\sigma_n(M+N-zI)\le\eps}+\int_{\abs z>R+\eps}\Pr(\abs z-\eps<\magn A) \wrt z
\\&\le\pi (R + \eps)^2\sup_{\abs z\le R+\eps}\Pr\pare{\sigma_n(M+N-zI)\le\eps}+2\pi\int_R^\infty (x + \eps)\Pr(x<\magn{M+N})\wrt x
\\&\le\pi (R + \eps)^2\sup_{\magn{M'}\le1.5R+\eps}\Pr\pare{\sigma_n(M'+N)\le\eps}+2\pi\int_{R-\magn M}^\infty (2x + \eps)\Pr(x<\magn{N})\wrt x.}
The final term is bounded by $\exp(-\Omega(R))$ since $\magn N$ has subgaussian concentration around $\sqrt n$, and $x$ is at least $R-\magn M\ge R/2\gg\sqrt n$.
\end{proof}

The following simple lemma, similar to net arguments appearing in \cite{b1,b3}, relates minimum gap bounds to singular value tail bounds.
\begin{lemma}[Spacing to smallest two singular values]\label{a25}
Let $A$ be any random matrix. 
Fix any constants $R>r>0$ and $\lambda\in(0,1)$.
Then
\spliteq{}{
\Pr\pare{\eta(A)\le r}&\le\frac{8(R+r)^2}{r^2}\sup_{\abs z\le R}\sqbrac{\Pr\pare{\sigma_n(A-zI)\le r^{1+\lambda}}+\Pr\pare{\sigma_{n-1}(A-zI)\le r^{1-\lambda}}} \\&+ \Pr(\specr\pare{M+N}\ge R).
}
\end{lemma}
\begin{proof}
Let 
\[
    \net=\left(\frac{r}{\sqrt 2}\Z+i\frac{r}{\sqrt 2}\Z\right)\cap\ball(0,R + r).
\]
Let $\lambda_i,\lambda_j$ be the eigenvalues of $A$ achieving $\abs{\lambda_i-\lambda_j}=\eta(A)$.
Consider the event that $\eta(A)\le r$.
If either of $\lambda_i,\lambda_j$ are outside $\ball(0,R)$, then by definition we have $\specr(A)\ge R$. On the other hand, if they are both inside the disk, then there must be some $z\in\net$ such that $\lambda_i,\lambda_j\in\ball(z,r)$.
So by union bound,
\[
\Pr\pare{\eta\pare{A}\le r}\le
\Pr\pare{\specr(A)\ge R}+\sum_{z\in\net}\Pr\pare{\abs{\lambda_i-z}\abs{\lambda_j-z}\le r^2}.
\]
Note that $\lambda_i-z$ and $\lambda_j-z$ are eigenvalues of $A-zI$, so by Weyl majorization (see \cite[Theorem 3.3.14]{b29}) we deterministically have
\spliteq{\label{a33}}{
\abs{\lambda_i-z}\abs{\lambda_j-z}\ge\sigma_n(A-zI)\sigma_{n-1}(A-zI).
}
Finally, the quantity (\ref{a33}) being upper bounded by $r^2$ implies one of the following events must occur:
\[\sigma_n(A-zI)\le r^{1+\lambda}\,\,\vee\,\,\sigma_{n-1}(A-zI)\le r^{1-\lambda}.\]
Therefore, for each $z\in\net$ we have by union bound
\spliteq{}{
\Pr\pare{\abs{\lambda_i-z}\abs{\lambda_j-z}\le r^2}
  &\le\Pr(\sigma_n(A-zI)\le r^{1+\lambda})+\Pr(\sigma_{n-1}(A-zI)\le r^{1-\lambda})
}
The result follows by noting $\abs\net\le 8 (R+r)^2/r^2$.
\end{proof}
We note that given bounds of type \eqref{a13}, the conclusion of Lemma \ref{a25} is nontrivial whenever
\begin{equation}\label{a34} \sup_{\lambda \in (0,1)}\min\pare{(1+\lambda)c_0,\, (1-\lambda)c_1} > 2.\end{equation}
The left hand side is optimized for $\lambda=(c_1-c_0)/(c_1+c_0)$ in which case \eqref{a34} reduces to \eqref{a14}.

\section{Bounds on singular values} \label{a8}
Section \ref{a28} showed that control over $\E \Lambda_\eps(A)$ and $\eta(A)$ can be reduced to control over the bottom two singular values $\sigma_n(z-A)$ and $\sigma_{n-1}(z-A)$.
As described in the introduction, this section shows tail bounds of the form \eqref{a13} for $m=0,1$ with $c_0=2$ and $c_1=4$. 

This is the only section where we strongly use the model of random matrix $A=M+N_g$. In particular, it is important that the entries of $N_g$ are independent and are of the form $\delta\cdot g$ where $\delta$ is Bernoulli and $g$ has an absolutely continuous distribution on $\C$. We introduce some helpful notation for this section. The L\'evy concentration (or small-ball probability) of a complex random variable $X$ is the function
\[\L(X, r) = \sup_{z \in \C}\Pr\pare{X\in\ball(z, r)}.\]
Let $\vec g$ and $\vec\delta_\rho$ denote a vector of i.i.d.~complex standard Gaussians and a vector of Bernoulli random variables with mean $\rho$, respectively.
For any complex vector $v$ we write
\[p(v, r) := \L(\rinr{\vec g \odot \vec \delta_\rho}{v}, r),\]
where $\odot$ is the Hadamard (entrywise) product.
This definition is motivated by the simple observation that the independence of the rows of $N_g$ implies
\[\Pr\pare{\magn{(M+N_g)v}_\infty\le r}\le p(v,r)^n.\]

We can partition the sphere according to value of $p(v,r)$ by defining
\[\alg{Incomp}(r,s)=\set{v\in\sph{n-1}:p(v,r)\le s},\quad\alg{Comp}(r,s)=\set{v\in\sph{n-1}:p(v,r) > s}.\]

The standard approach for establishing a singular value bound is to show by casework that
it is unlikely for any vector in $\alg{Incomp}$ or in $\alg{Comp}$ to witness the small singular space.

We begin with the $\alg{Incomp}$ case. The first fact establishes a basis for the small singular space of a certain form. The strong anti-concentration properties of $\alg{Incomp}$ then yield the result.

\begin{fact}[Reduced row echelon form basis]\label{a19}
Every dimension $k$ subspace of $\mathbb C^n$ can be expressed as the column space of a matrix of the form 
\[
    P\bmat{D\\X},
\]
where $P \in \mathbb C^{n \times n}$ is a permutation matrix, $D \in \mathbb C^{k \times k}$ is diagonal with $|D_{jj}| \ge 1/\sqrt n$ for all $j \in [k]$, and $X \in \mathbb C^{(n-k) \times k}$ such that the columns of
$\left[\begin{smallmatrix} D \\ X \end{smallmatrix}\right]$
are unit length.
\end{fact}
This fact follows from running Gaussian elimination with complete pivoting on a basis and normalizing. We now have

\begin{lemma}[Incompressible vectors do not witness small singular space]\label{a20}
Fix any $\eps,s\in(0,1)$ and $m \in \Z_{\ge0}$ with $m < n$. Then
    \[
    \Pr\pare{\sigma_{n-m}(M+N_g)\le\eps\wedge\exists v\in\alg{Incomp}(2\eps\sqrt n,s)\,\st\,\magn{(M+N_g)v}\le\eps}\le
    \binom{n}{m+1}s^{m+1}.
    \]
\end{lemma}
\begin{proof}
For the moment suppose $\sigma_{n-m}(M+N_g)\le\eps$. Then there exists a subspace $W$ of dimension $m+1$ such that $\magn{w^*(M+N_g)}\le\eps$ for each $w\in W\cap\sph{n-1}$. By Fact \ref{a19}, there exists a permutation $P$, diagonal $D$ with $\abs{D_{jj}}\ge1/\sqrt n$, and some $X$ such that
\begin{equation} \label{a35}
    \max_j\magn{e_j^*\bmat{D\\X}^*P^*(M+N_g)}\le\eps.
\end{equation}

Let $\mathcal E$ denote the event
\[
    \sigma_{n-m}(M+N_g)\le\eps\wedge\exists v\in\alg{Incomp}(2\eps\sqrt n,s)\,\st\,\magn{(M+N_g)v}\le\eps.
\]
For permutation matrix $P$ let $\mathcal E_P$ denote the event \eqref{a35} for some $D, X$ satisfying the aforementioned properties. Let $\mathcal U$ be a function that takes $B \in \mathbb C^{(n-m-1)\times n}$ and chooses a unit vector $u$ such that $p(u, 2\eps \sqrt{n}) \le s$ and $\|Bu\| \le \eps$ (and is undefined if no such vector exists).

Now fix some permutation matrix $P$. We can write $P^*(M+N_g)=\left[\begin{smallmatrix}B_1 \\ B_2\end{smallmatrix}\right]$ where $B_1$ has $m+1$ rows and $B_2$ has $n-m-1$ rows. Observe $B_1$ and $B_2$ are independent random matrices. Furthermore, if $\mathcal E$ occurs then there exists unit vector $v$ such that
\[\magn{(M+N_g)v}=\sqrt{\magn{B_1v}^2+\magn{B_2v}^2}\le\eps\]
and $p(v, 2\eps \sqrt{n}) \le s$, so $\mathcal U$ is well-defined on $B_2$.
Hence let $u = \mathcal U(B_2)$ and observe the event $\mathcal E \wedge \mathcal E_P$ implies
\spliteq{}{
\eps
   \ge\max_j\magn{e_j^*\bmat{D\\X}^*P^*(M+N_g)}
  &=  \max_j\magn{e_j^*(D^*B_1+X^*B_2)}
\\&\ge\max_j\magn{e_j^*D^*B_1u+e_j^*X^*B_2u}
\\&\ge\max_j\magn{D_{jj}^*e_j^*B_1u+e_j^*X^*B_2u}.
}
Because $\magn{B_2u}\le\eps,\magn{e_j^*X}\le1,\abs{D_{jj}}\ge1/\sqrt{n}$, we conclude
\spliteq{}{
2\eps\sqrt n\ge\max_j\magn{e_jB_1u}=\magn{B_1u}_\infty.
}
Consequently, we have
\[
    \Pr\pare{\mathcal E \wedge \mathcal E_P} = \mathbb E \left[ \Pr \left[\mathcal E \wedge \mathcal E_P | B_2 \right] \right] \le \mathbb E\left[ \Pr\pare{\magn{B_1u}_\infty \le 2\eps \sqrt{n}} | B_2\right] = \mathbb E\left[p(u,2\eps\sqrt n)^{m+1} | B_2 \right] \le s^{m+1},
\]
where the second equality follows because $u$ depends only on $B_2$ and hence is chosen irrespective of $B_1$. To bound $\Pr\pare{\mathcal E}$ it remains to sum over permutations $P$. However, despite fixing $P$ earlier the only information we used about $P$ was which indices correspond to blocks $B_1$ and $B_2$. Hence it is only necessary to pay a union bound cost of $\binom{n}{m+1}$, yielding the desired bound.
\end{proof}

For the $\alg{Comp}$ case, we no longer have strong anti-concentration of $\rinr{\vec g\odot\vec\delta_\rho}{v}$. Instead we show $\alg{Comp}$ has small entropy and and use a net argument.

\begin{lemma}[Compressible vectors are close to sparse]\label{a21}
If $p(v,r) > s$, then
\[\frac{\log(2/s)}\rho\ge\#\set{j:\abs{v_j}\ge r\sqrt{2/s}}.\]
\end{lemma}
\begin{proof}
By iterated expectation we have
\spliteq{}{
s < p(v,r)
  &=\Pr\pare{\abs{\inr{\vec g\odot\vec\delta}{v}}\le r}
\\&=\E_\delta\Pr_g\pare{\abs{\inr{\vec g}{\vec\delta\odot v}}\le r}
\\&\le\E_\delta\min\pare{1,\frac{r^2}{\rmagn{\vec\delta\odot v}^2}}.
\\&\le\Pr\pare{\rmagn{\vec\delta\odot v}^2\le\frac{2r^2}s}+\frac s2.
\\&\le(1-\rho)^{\#\set{j:\abs{v_j}^2\ge 2r^2/s}}+\frac s2.
}
Using $1-\rho\le e^{-\rho}$ and rearranging gives the result.
\end{proof}

\begin{lemma}[Compressible vectors do not witness small singular space]\label{a22}
Suppose $\frac{1}{n} < \rho \le 1$, $0 < s < \frac{2}{e}$, $\eps > 0$, and integer $K \ge 1$. Let
\[
    T = \left(2 \eps n \sqrt{\frac{2}{s}}\right)^{1/K}, \quad H = \left(\frac{\log(2/s)}{\rho}\right)^{1/K}.
\]
Furthermore let $R \ge 1$ and suppose
\begin{equation} \label{a36}
    T \le \frac{1}{Rn}, \quad \frac{\eps}{T^{(K-1)}} \le \frac{1}{n^2}, \quad T \ge (4n)^{1/K}\exp\left(-\frac{\rho n }{4KH}\right), \quad H^K \le n, \quad H^{K-1} \le \frac{1}{\rho}.
\end{equation}
Then for $n \ge 30$ we have
\begin{equation} \label{a37}
    \Pr\pare{\exists v\in\alg{Comp}(2\eps\sqrt n,s)\,\st\,\magn{(M+N_g)v}\le\eps}
    \le 55K^2 e^{-\frac34 n\rho} +\Pr\pare{\magn{M+N_g}_\infty\ge R}.
\end{equation}
\end{lemma}
\begin{proof}
Suppose such $v$ exists.
We consider a range of thresholds from $1/\sqrt{n}$ decreasing to $2\eps\sqrt n\cdot\sqrt{2/s}$ and estimate how many entries of $v$ are above each threshold.
The key parameter is the number of thresholds we consider. In this proof we consider $K+1$ exponentially spaced thresholds given by
\[
    \frac{T^j} {\sqrt n}, \quad 0 \le j \le K.
\]
These thresholds are decreasing in $j$ because $T \le 1/(Rn) < 1$.
Let $a_j(v)$ be the number of entries of $v$ whose magnitude is above the $j$th threshold, i.e.
\[
    a_j(v)=\#\set{i:\abs{v_i}\ge\frac1{\sqrt n} T^j}.
\]
Lemma \ref{a21} immediately implies that $a_K\le\frac{\log(2/s)}\rho$. On the other hand, $a_0(v)\ge1$ since $v$ is a unit vector.
So\[\frac{a_0(v)}{a_1(v)}\cdot\frac{a_1(v)}{a_2(v)}\cdots\frac{a_{K-1}(v)}{a_K(v)}=\frac{a_0(v)}{a_K(v)}\ge\frac\rho{\log(2/s)} = H^K.\]
Note $H > 1$. Therefore by the pigeonhole principle there exists indices $0 \le k <K$ and $0 \le j < K$ such that
\begin{equation} \label{a38}
    H^j \le a_{k}(v) \le a_{k+1}(v) \le H^{j+1}.
\end{equation}
Let us now fix $k, j$ and consider the set of vectors $v$ satisfying both \eqref{a37} and \eqref{a38}.
Given $v$, define $\tilde v$ be the vector obtained by first rounding the entries of $v$ whose magnitude are less than $ T^{k+1} / \sqrt{n}$ down to 0 and then rounding the real and imaginary parts of all other entries \textit{up} to the nearest multiple of $T^{k+1} / (2 \sqrt{n})$.
Denote the set of all possible $\tilde v$ as $\net_{k,j}$.
Then
\[
    \magn{v-\tilde v}_\infty\le\frac1{\sqrt n} T^{k+1}, \quad
    \quad a_{k}(v) \le a_{k}(\tilde v) \le a_{k+1}(\tilde v) = a_{k+1}(v).
\]
Thus
\spliteq{}{
\label{a39}
\magn{(M+N_g)v}\le\eps
  &\implies\magn{(M+N_g)v}_\infty\le\eps
\\&\implies\magn{(M+N_g)\tilde v}_\infty\le\eps+\magn{M+N_g}_\infty\magn{v-\tilde v}_\infty
\\&\implies\magn{(M+N_g)\tilde v}_\infty\le\eps+R T^{k+1}/\sqrt{n} \ \text{ or } \ \magn{M+N_g}_\infty\ge R.
}
Moreover for $\tilde v\in\net_{k, j}$ we have
    \spliteq{}{
    \Pr\pare{\magn{(M+N_g)\tilde v}_\infty\le\eps+\frac R{\sqrt n}T^{k+1}}
      &\le p\pare{\tilde v,\eps+\frac R{\sqrt n}T^{k+1}}^n
    \\&\le\sqbrac{(1-\rho)^{H^j}+\pare{\frac{\eps+\frac R{\sqrt n}T^{k+1}}{\frac{1}{\sqrt{n}}T^k}}^2}^n
    \\&=\sqbrac{(1-\rho)^{H^j}+\pare{\sqrt{n}\frac{\eps}{T^k}+ RT}^2}^n
    \\&=\sqbrac{(1-\rho)^{H^j}+\pare{\sqrt{n}\frac{\eps}{T^{K-1}}+ RT}^2}^n
    \\&\le\sqbrac{(1-\rho)^{H^j}+\frac{4}{n^2}}^n
    }
where the last step follows from the first and second conditions in \eqref{a36}. We now note that if $(1 - \rho)^{H^j} > 1/n$ then
\[
    \sqbrac{(1-\rho)^{H^j}+\frac{4}{n}}^n \le \sqbrac{(1-\rho)^{H^j}\left(1 + \frac{4}{n}\right)}^n \le e^4 (1 - \rho)^{n H^j}.
\]
On the other hand, if $(1 - \rho)^{H^j} \le 1/n$ then we have the upper bound $(5/n)^n$. Hence
\begin{equation}
    \label{a40}
    \sqbrac{(1-\rho)^{H^j}+\frac{4}{n}}^n \le
    \min \left(e^4 (1 - \rho)^{n H^j}, \left[\frac{5}{n}\right]^n\right)
    \le
    \min \left(e^{-\rho n H^j + 4}, \left[\frac{5}{n}\right]^n\right).
\end{equation}

We now bound the cardinality of $\net_{k,j}$.
Because $a_k(\tilde v) \le H^{j+1} \le n$ we have
\begin{equation}
    \label{a41}
    \begin{split}
        \abs{\net_{k, j}}
          &\le\binom{n}{H^{j+1}}\pare{\frac{4\sqrt n}{T^{k+1}}}^{2 H^{j+1}}\\
          &\le\binom{n}{H^{j+1}}\pare{\frac{4\sqrt n}{T^{K}}}^{2 H^{j+1}}\\
          &\le\binom{n}{H^{j+1}}\pare{\frac{16 n}{T^{2K}}}^{H^{j+1}}\\
          &\le\pare{\frac{16 n^2}{T^{2K}}}^{H^{j+1}}\\
          &\le\exp\pare{\frac{\rho n H^j}{4}}
    \end{split}
\end{equation}
where the last inequality follows from the third condition in \eqref{a36}.
By combining \eqref{a41} with \eqref{a39}, \eqref{a40}, then applying the union bound first over $\tilde v \in \net_{k,j}$ and then over the indices $k,j$ we conclude
\begin{align*}
    &\Pr\pare{\exists v\in\alg{Comp}(2\eps\sqrt n,s)\,\st\,\magn{(M+N_g)v}\le\eps}\\
    &\quad\quad\quad\quad\quad\quad\le \sum_{k,j} \exp\pare{\frac{\rho n H^j}{4}} \min \left(e^{-\rho n H^j + 4}, \left[\frac{5}{n}\right]^n\right) + \bP\left(\|M+N_g\|_\infty \ge R\right)\\
    &\quad\quad\quad\quad\quad\quad\le \sum_{k,j} e^{-3\rho n H^j/4 + 4} + \sum_{k,j} e^{\rho n H^{K-1}}
    \left[\frac{5}{n}\right]^n + \bP\pare{\|M+N_g\|_\infty \ge R}\\
    \intertext{Using the fifth condition in \eqref{a36} and the condition $n > 30$ we obtain}
    &\quad\quad\quad\quad\quad\quad\le \sum_{k,j} e^{-3\rho n H^j/4 + 4} + \sum_{k,j} \left[\frac{5e}{n}\right]^n + \bP\left(\|M+N_g\|_\infty \ge R\right)\\
    &\quad\quad\quad\quad\quad\quad\le K^2 \exp\pare{-\frac{3\rho n }{4} + 4} + K^2 \left[\frac{5e}{n}\right]^n + \bP\left(\|M+N_g\|_\infty \ge R\right)\\
    &\quad\quad\quad\quad\quad\quad\le 55 K^2 \exp\pare{-\frac{3\rho n }{4}} + \bP\left(\|M+N_g\|_\infty \ge R\right).
\end{align*}
\end{proof}
To apply Lemma \ref{a22} we need to make a prudent choice of the parameters $\eps, s, \rho, K$, and $R$ so that \eqref{a36} is satisfied. One such choice is given by
\begin{lemma} \label{a42}
    Fix $R \ge 4$ and integer $K \ge 1$. The conditions of Lemma \ref{a22} are satisfied when
    \[
        \frac{1}{R^{20(K+1)} n^{20(K+1)}} \le s < \frac{2}{e}, \quad \eps^2 = \frac{s}{8R^{2K}n^{2(K+1)}},
    \]
    \[
        \frac{n^{1/(K+1)}}{n} 20 (K+2) \log (Rn) \le \rho \le \left(\frac{1}{20 (K+2)\log (Rn) }\right)^{(K-1)}
    \]
\end{lemma}
\begin{proof}
    It is easy to see the trivial conditions $1/n < \rho \le 1$, $0 < s < 2/e$, $\eps > 0$ so it remains to check \eqref{a36}. We note
    \begin{equation} \label{a43}
        T = \left(2 \eps n \sqrt{\frac{2}{s}}\right)^{1/K} = \frac{1}{Rn},
    \end{equation}
    which satisfies the first condition in \eqref{a36}. Moreover, \eqref{a43} implies
    \[
        \frac{1}{n^4} T^{2(K-1)} = \frac{1}{R^{2(K-1)}n^{2(K+1)}} \ge \frac{1}{4e R^{2K}n^{2(K+1)}} \ge \eps^2,
    \]
    which satisfies the second condition in \eqref{a36}.
    Now observe
    \begin{equation} \label{a44}
        \rho n = \left( \frac{\rho^{K+1}n^K}{\rho} \right)^{1/K}
        = H \left( \frac{\rho^{K+1}n^K}{\log(2/s)} \right)^{1/K}
        \ge H \left( \frac{\rho^{K+1}n^K}{20(K+2)\log(Rn)} \right)^{1/K}
        \ge 20 H(K+2)\log(Rn).
    \end{equation}
    Hence
    \[
        T = \frac{1}{Rn} = (4n)^{1/K}\exp\left( -\frac{4K H\log(Rn) + H \log(4n)}{4KH}\right) \ge (4n)^{1/K}\exp\left( -\frac{\rho n}{4KH}\right),
    \]
    fulfilling the third condition. The fourth condition is satisfied because
    \begin{equation} \label{a45}
        H^{K-1} \le \left(\frac{20(K+2)\log(Rn)}{\rho}\right)^{(K-1)/K} = \frac{1}{\rho}\left(\rho^{1/(K-1)}20(K+2)\log(Rn)\right)^{(K-1)/K} \le \frac{1}{\rho}.
    \end{equation}
    Finally we observe \eqref{a44} and \eqref{a45} imply  $H^{K} \le H/\rho \le n \rho / \rho \le n$.
\end{proof}

Combining Lemma \ref{a20}, Lemma \ref{a22}, and Lemma \ref{a42} then yields
\begin{proposition} \label{a23}
Let $n$ be sufficiently large. Choose $R, \rho, \eps$ such that
\[
    \frac{(\log R)^{3}(\log n)^{3}}{n} \le \rho \le \frac{n^{1/2}}{n}, \quad K = \left\lceil \frac{3\log(n)}{2\log(n\rho)} \right\rceil,
    \quad
    \frac{1}{n^{10(K+1)}R^{10(K+1)}} \le \eps\le\frac{1}{4 n^{K+1}R^{K}}.
\]
Then
\[
\Pr
\pare{\sigma_{n-m}(M+N_g)\le\eps}\le \left[ 8 n^{2K+3}R^{2K}\eps^{2} \right]^{m+1} + n^2 \exp\pare{-\frac{3}{4}\rho n}+\Pr\pare{\magn{M+N_g}_\infty\ge R}.
\]
\end{proposition}
\begin{proof}
    Set
    \[
        s = 8R^{2K}n^{2(K+1)}\eps^2.
    \]
    Assume the parameters $R, \eps, s, R, \rho$ satisfy the conditions in Lemma \ref{a42}. Then Lemma \ref{a22} applies.
    Because $\alg{Comp}(2\eps\sqrt n,s)\sqcup\alg{Incomp}(2\eps\sqrt n,s)$ is a partitioning the sphere, the event $\sigma_{n-m}(M+N_g)\le\eps$ is the union of the events whose probabilities are bounded in each of Lemma \ref{a20} and Lemma \ref{a22}. We obtain the final result by a union bound.

    It remains to check the conditions in Lemma \ref{a42}. The first two conditions are immediate and the last condition follows from
    \begin{align*}
        \frac{n^{1/(K+1)}}{n} 20 (K+2) \log (Rn)
        &\le \frac{n^{2\log(n\rho)/(3\log n)}}{n}12(K+2)\log(R)\\
        &\le \frac{(n\rho)^{2/3}}{n}40(K+2)\log(R)\\
        &\le \rho \frac{1}{(n\rho)^{1/3}}40(K+2)\log(R)\\
        &\le \rho \frac{40(K+2)}{\log n} \le \rho
        \intertext{and}
        \left(\frac{1}{20(K+2)\log (Rn) }\right)^{K-1}
        &\ge \left(\frac{1}{40(K+2)\log (R) }\right)^{K-1}\\
        &\ge \left(\frac{1}{50\log n \log (R) }\right)^{K-1}\\
        &\ge \left(\frac{1}{50\log n \log (R) }\right)^{3\log n / (6\log \log n + 6\log \log R)}\\
        &\ge \exp\left(\log 50 -\frac{\log n}{\log \log n}\right) n^{-1/2} \ge n^{-1/2} \ge \rho.
    \end{align*}
\end{proof}
\begin{remark}
For $m>0$, there is overlap between the two events we union over in Proposition \ref{a23} and therefore slack in the bound. In particular, when $\sigma_{n-m}(M+N_g)$ is small, there should be an entire $(m+1)$-dimensional space witnessing the small singular space. Lemma \ref{a20} bounds the probability that $\sigma_{n-m}(M+N_g)$ is small and \textit{a single} vector in $\alg{Incomp}$ witnesses the small singular space.
Lemma \ref{a22} bounds the probability \textit{a single} vector in $\alg{Comp}$ witnesses the small singular space.
But it is possible that the $(m+1)$-dimensional small singular space contains both vectors in $\alg{Incomp}$ and $\alg{Comp}$.
\end{remark}
\section{Proof of Theorem \ref{a1}}
We now prove the main theorem of this paper.

\begin{proof}[Proof of Theorem \ref{a1}]
Let
\[
    f_{R'}(m,\eps)=\sup_{\magn{M'}\le R'}\Pr(\sigma_{n-m}(M'+N)\le\eps).
\]
Let $R = 2\magn M+n$ so that $\Pr(\magn{M+N}\le R)\le e^{-\Omega(n)}$.
Then combining the reductions Lemma \ref{a24} and \ref{a25} with Proposition \ref{a18} gives
\allowdisplaybreaks
\begin{align*}
    \Pr\pare{\max_j \kappa(\lambda_j)^2\ge\tau} &\le 8(R +\eta)^2\pare{\frac{f_R(0,\eta^{4/3})+f_R(1,\eta^{2/3})}{\eta^2}}\\
    &\quad\quad\quad\quad+\left[3n\lg\pare{\frac{3R}{\eta}} + 3\lg\pare{\frac{2\eta}{\tau}} \right]\frac{n^2f_{R/2+\eta/(n\sqrt{\tau})}\pare{0,\frac\eta{n\sqrt\tau}}}{\eta^2}+e^{-\Omega(n)}.
\end{align*}
for any $\tau > 0$ and $R > \eta > 0$.
Note the first term on the right hand side is an upper bound for $\Pr(\eta(M+N) \le \eta)$.
Set $\tau=\eta^{-2/3}/n^2$. Then we have
\begin{align*}
    \Pr\pare{\max_j \kappa(\lambda_j)^2 \ge \frac1{n^2\eta^{2/3}}}
    &\le 8(R +\eta)^2\pare{\frac{f_R(0,\eta^{4/3})+f_R(1,\eta^{2/3})}{\eta^2}}\\
    &\quad\quad\quad\quad+\left[3n\lg\pare{\frac{3R}{\eta}} + 3\lg\pare{2\eta^{5/3}n^2} \right] \frac{n^2f_{2R+\eta^{4/3}}\pare{0,\eta^{4/3}}}{\eta^2}+e^{-\Omega(n)}.
\end{align*}
Now we would like to apply Proposition \ref{a23} with $m=0$ and $m=1$ to bound $f(m,\eps)$.
To do so we set
\[
    \eta^2 = \frac{1}{n^{15(K+1)}R^{15(K+1)}}.
\]
Then Proposition \ref{a23} implies
\begin{align*}
    \frac{f_R(0, \eta^{4/3})}{\eta^2}
    &\le 8n^{2K + 3}R^{2K}\eta^{2/3} + \frac{1}{\eta^2}\left[n^2\exp\left(-\frac{3}{4}\rho n\right) + e^{-\Omega(n)}\right]\\
    &\le n^{-5K}R^{-5K} + \frac{n^2}{\eta^2}\exp\left(-\frac{3}{4}\rho n\right) + e^{-\Omega(n)} \le n^{-2K}R^{-2K} + e^{-\Omega(\rho n)}\\
    \frac{f_R(1, \eta^{2/3})}{\eta^2}
    &\le 64n^{4K + 6}R^{4K}\eta^{2/3} + \frac{1}{\eta^2}\left[n^2\exp\left(-\frac{3}{4}\rho n\right) + e^{-\Omega(n)}\right] \le 64n^{-K + 1}R^{-K}+ e^{-\Omega(\rho n)}\\
    \frac{f_{R/2 + \eta^{4/3}}(0, \eta^{4/3})}{\eta^2}
    &\le
    \frac{f_{R}(0, \eta^{4/3})}{\eta^2} \le
     n^{-2K}R^{-2K} + e^{-\Omega(\rho n)}.
\end{align*}
Hence we conclude
\[
    \Pr\pare{\max_j \kappa(\lambda_j)^2} \le O\left(\frac{1}{n^{2K}}\right) + e^{-\Omega(\rho n)}.
\]
Note the first term dominates. Similarly we have
\spliteq{}{
\Pr\pare{\eta(M+N) \le \eta}
\le
8(R +\eta)^2\pare{\frac{f_R(0,\eta^{4/3})+f_R(1,\eta^{2/3})}{\eta^2}}
\le O\left(\frac{1}{n^{2K}}\right) + e^{-\Omega(\rho n)}.}

\end{proof}

\section{Application to GMRES}
\label{a46}
As alluded to in the introduction, controlling the eigenvector condition number enables the analysis of algorithms for non-Hermitian inputs. It is particularly important for eigenvalue problems, and yet our application shows its importance to solving systems of equations. In particular GMRES constructs polynomial approximations to the inverse function, so reduces system-solving to a problem relating to the spectrum of the input. The work of \cite{b2} achieves regularity of the input by adding a dense matrix $N$ to the input. 
The drawback is that the cost of computing matrix-vector products goes from $O(\nnz(M))$ to $O(\nnz(M)+\nnz(N))=\Theta(n^2)$ where $\nnz(A)$ is the number of nonzero entries in the matrix $A$. The algorithmic content of this paper is that it suffices to take $N$ to be a sparse perturbation, with $\E\nnz(N)=n^2\rho$ for $\rho=\Theta(\log^6(n)/n)$.

\begin{theorem}\label{a10}
Let $\Omega \subset \mathbb C$ be a set such that
\[\inf_{\deg(p)\le k}\sup_{z\in\Omega}\abs{1-zp(z)}\le C\alpha^k.\]
Let $\eps, \delta \in (0,1)$ and let
$M$ be a matrix such that $\magn M\le\poly(n)$ and
$\Lambda_\delta(M)\subset\Omega$. \alg{shatteredGMRES} is an algorithm which, in \[k=O\pare{\dfrac{\log(n/\delta)\log(n)+\log(1/\eps)}{\log(1/\alpha)}}\] iterations, outputs $x$ such that
\[(M+\Delta_M)x=(b+\delta_b), \quad \magn{\Delta_M}\le\delta,\quad\magn{\delta_b}\le\eps\]
and uses
\(O\pare{ k\nnz(M)+kn\log^6(n)+k^2n }\)
arithmetic operations.
\end{theorem}
\begin{proof}
Sample $N_g$ with sparsity $\rho=\log(n)^6/n$ and set $A=M+(0.1\delta/{\sqrt n})N_g$.
\alg{shatteredGMRES} is GMRES applied to $A$.
By the minimum residual property, \cite[Proposition 4]{b10}, the output $x$ of GMRES applied to $A$ with $k$ iterations satisfies
\begin{equation}\label{a47}
\frac{\magn{Ax-b}}{\magn b}
\le\kappa_V(A)\inf_{\deg(p)\le k}\max_{j=1,\ldots,n}\abs{1-\lambda_jp(\lambda_j)},\end{equation}
where $\lambda_j$ are the eigenvalues of $A$.
Note $\lambda_j\in\Lambda_{0.1\frac{\delta}{\sqrt{n}}\magn{N_g}}$, and by concentration of $\magn{N_g}$, $\Lambda_{0.1\frac\delta{\sqrt n}\magn{N_g}}\subset \Lambda_\delta(M)$ with high probability, and by hypothesis $\Lambda_\delta(M)\subset\Omega$.
As a consequence, \eqref{a47} becomes
\[\frac{\magn{Ax-b}}{\magn b}\le C\kappa_V(A)\alpha^k\]
Control over $\kappa_V(A)$ is implied by Corollary \ref{a7} as $\log\kappa_V(A)\le O\pare{\log(n/\delta)\log(n)}$ with high probability.
Notice the average cost of each iteration is
\[O\pare{\nnz(A) + kn}
=O\pare{\nnz(M)+n^2\rho + kn}\]
giving the final result.
\end{proof}
A convenient application of this theorem is to $\Omega=\ball(z_0,r)$ for $\abs{z_0}>r$, in which case $\alpha=r/\abs{z_0}$.

\bibliographystyle{alpha}
\bibliography{outbib}
\end{document}